\documentclass[12pt]{amsart}
\usepackage{amssymb,amsmath,amscd,graphicx,
latexsym,amsthm}
\usepackage{amssymb,latexsym,eufrak,amsmath,amscd,graphicx}
  \usepackage[all]{xy}
\theoremstyle{plain}
\newtheorem{theorem}{Theorem}[section]
\newtheorem{classicaltheorem}[theorem]{Classical
Theorem}
\newtheorem{proposition}[theorem]{Proposition}
\newtheorem{corollary}[theorem]{Corollary}

\theoremstyle{definition}
\newtheorem{definition}[theorem]{Definition}

\newtheorem{remark}[theorem]{Remark}
\newtheorem{example}[theorem]{Example}

\newtheorem{problem}[theorem]{Problem}
\newtheorem{conjecture}[theorem]{Conjecture}

\newcommand{\lra}{\longrightarrow}
\newcommand{\noi}{\noindent}
\newcommand{\PP}{\mathbf{P}}

\newcommand{\RR}{\mathbf{R}}

\newcommand{\NN}{\mathbf{N}}
\newcommand{\ZZ}{\mathbf{Z}}
\newcommand{\CC}{\mathbf{C}}
\newcommand{\QQ}{\mathbf{Q}}

\newcommand{\OO}{\mathcal{O}}

\newcommand{\fra}{\mathfrak{a}}
\newcommand{\frb}{\mathfrak{b}}

\newcommand{\bull}{_{\bullet}}

\newcommand{\eps}{\varepsilon}

\newcommand{\lbl}{\vskip 6pt}

\newcommand{\SBL}[1]{\mathbf{B}({#1})}
\newcommand{\BAug}[1]{\mathbf{B}_{+}({#1})}

\newcommand{\HH}[3]{H^{{#1}} \big( {#2} , {#3}
\big) }
\newcommand{\hh}[3]{h^{{#1}} \big( {#2} , {#3}
\big) }

\newcommand{\Effbar}{\overline{\textnormal{Eff}}}
\newcommand{\Nef}{\textnormal{Nef}}

\newcommand{\ord}{\textnormal{ord}}

\newcommand{\pr}{\prime}

\newcommand{\Div}{\text{Div}}

\newcommand{\vol}{\textnormal{vol}}

\newcommand{\Bl}{\text{Bl}}

\newcommand{\Linser}[1]{| \mspace{1.5mu} {#1}
\mspace{1.5mu} |}
\newcommand{\linser}[1]{\Linser{  {#1}  }}
\newcommand{\alinser}[1]{\| \mspace{1.5mu}{#1}
		\mspace{1.5mu}\|}
\newcommand{\Bs}[1]{\textnormal{Bs}\big(
\mspace{1.1mu} | {#1} | \mspace{1.1mu} \big)}
\newcommand{\bs}[1]{\frb \big(\mspace{1.1mu}
\linser{{#1}} \mspace{1.1mu} \big)}

\newcommand{\hhat}[2]{\widehat{h}^{{#1}}\big(
{#2} \big)}

\begin{document}

\title{Asymptotic invariants of line bundles} 

\author[L.Ein]{Lawrence Ein}
\address{Department of Mathematics \\ University
of Illinois at Chicago, \hfil\break\indent  851
South Morgan Street (M/C 249)\\ Chicago, IL
60607-7045, USA}
\email{ein@math.uic.edu}

\author[R. Lazarsfeld]{Robert Lazarsfeld} 
\address{Department of Mathematics
\\ University of Michigan \\ Ann Arbor, MI
48109, USA}
\email{rlaz@umich.edu}

\author[M. Musta\c{t}\v{a}]{Mircea
Musta\c{t}\v{a}}
\address{Department of Mathematics \\
University of Michigan \\   Ann Arbor, MI
48109,  USA}
\email{mmustata@umich.edu}

\author[M. Nakamaye]{Michael Nakamaye}
\address{Department of Mathematics and
Statistics\\ University of New
  Mexico,
\hfil\break\indent Albuquerque New Mexico 87131,
USA}
\email{nakamaye@math.unm.edu}

\author[M. Popa]{Mihnea Popa}
\address{Department of Mathematics \\
  Harvard University \\ 1 Oxford Street,
Cambridge,
  \hfil\break\indent MA 02138,
  USA}
\email{mpopa@math.harvard.edu}

\thanks{The research of the authors was
partially supported by the NSF under grants
 DMS 
0200278,   DMS
  0139713,   DMS 0500127, DMS
  0070190, and   DMS 0200150.}

\maketitle

\section*{Introduction}

Let $X$ be a smooth complex projective variety
of dimension $d$. It is classical that ample
divisors on $X$ satisfy many beautiful
geometric, cohomological, and numerical
properties that render their behavior
particularly tractable. By contrast, examples
due to Cutkosky and others (\cite{Cutkosky86a},
\cite{Cutkosky-Srinivas93b},
\cite[Chapter 2.3]{PAG}) have led to the
common impression that the linear series
associated to non-ample effective divisors are
in general   mired in pathology.

However, starting with fundamental work of
Fujita
\cite{Fujita94a}, Nakayama \cite{Nakayama04},
and Tsuji
\cite{Tsuji99}, it has recently become apparent
that arbitrary effective (or ``big") divisors  
in fact display a surprising number of
properties analogous to those of ample line
bundles.\footnote{The published record gives a
somewhat misleading sense of the chronology
here. An early version of Nakayama's  2004
memoir
\cite{Nakayama04} has been circulating as a
preprint since 1997, and Tsuji's ideas involving
asymptotic invariants occur in passing in a
number of so far unpublished preprints dating
from around 1999. Similarly, the results from
\cite{PAG} on the function
$\vol_X$ that we discuss below  initially
appeared in a preliminary draft of
\cite{PAG} circulated in 2001.} The key is to
study the properties in question from an
asymptotic perspective. At the same time, many
interesting questions and open problems
remain.  

The purpose of the present expository
note is to give an invitation to this
circle of ideas.  Our hope is that this informal
overview  might serve as a jumping off
point for  the more technical literature
in the area. Accordingly, we sketch many
examples but include no proofs.    In an
attempt to make the story as appealing as
possible to non-specialists, we focus on one
particular invariant --- the ``volume" --- that
measures the rate of growth of sections of
powers of a line bundle. Unfortunately, we must
then content  ourselves with giving references
for a considerable amount of   related work. The
papers
  \cite{Boucksom02},
\cite{Boucksom04}  of Boucksom from the
analytic viewpoint, and the exciting results of
Boucksom--Demailly--Paun--Peternell \cite{BDPP}
deserve particular mention: the reader can
consult 
\cite{Debarre05} for a survey.

We close this introduction by recalling some
notation and basic facts about cones of
divisors. Denote by $N^1(X)$
the N\'eron--Severi group of numerical
equivalence classes of divisors on $X$, and by
$N^1(X)_{\QQ}$ and $N^1(X)_{\RR}$ the
corresponding finite-dimensional rational and
real vector spaces  parametrizing numerical
equivalence classes of $\QQ$- and
$\RR$-divisors respectively. The N\'eron--Severi
space $N^1(X)_{\RR}$ contains two important
closed convex cones:
\begin{equation}
  N^1(X)_{\RR} \ \supseteq \
\Effbar(X)
\
\supseteq \ \Nef(X). \tag{*} \end{equation}
The
\textit{pseudoeffective cone}
$\Effbar(X)$
is defined to be the closure of the convex cone
spanned by the classes of all effective
 divisors on $X$. The \textit{nef  cone}
$\Nef(X) $ is the set of all nef divisor
classes, i.e. classes
$\xi \in N^1(X)_{\RR}$ such that $\big ( \xi
\cdot C \big) \ge 0$ for all irreducible curves
$C \subseteq X$. A very basic fact --- which in
particular explains the inclusions in (*)  ---
is that
\begin{align*}  \textnormal{interior} \Big(
\Effbar(X)
\Big)
\ &= \
\textnormal{Big(X)} \\
 \textnormal{interior} \Big( \Nef(X) \Big) \ &=
\ \textnormal{Amp(X)}. \end{align*}
Here $\textnormal{ Amp(X)}$ denotes the set of
ample classes on $X$, while
$\textnormal{Big(X)}$ is the cone of big
classes.\footnote{Recall that an integral
divisor $D$ is \textit{big} if
$\hh{0}{X}{\OO_X(mD)}$ grows like $m^{\dim X}$
for $m \gg 0$. This definition extends in a
natural way to $\QQ$-divisors (and with a
little more work to $\RR$-divisors), and one
shows that bigness depends only on the
numerical class of a divisor.} We refer to
\cite[Sections 1.4.C, 2.2.B]{PAG} for a fuller
discussions of these definitions and results.
Using this language, we may say that the theme
of this note is to understand to what extent
some classical facts about ample classes extend
to $\textnormal{Big}(X)$.

We are grateful to Michel Brion, Tommaso De
Fernex and Alex K\"uronya for valuable
discussions and suggestions.

\section{Volume of a line bundle}

In this section we give the definition of the
volume of a divisor and discuss its meaning in
the classical case of ample divisors. As before
$X$ is a smooth complex projective variety of
dimension
$d$.\footnote{For the most part, the
non-singularity hypothesis on $X$ is
extraneous. We include it here  
 only in order to avoid occasional
technicalities.}

Let   $D$ be a divisor on $X$. The
invariant on which we'll focus originates in
the 
\text{Riemann--Roch
problem} on $X$, which   asks for the
computation of the dimensions 
\[  \hh{0}{X}{\OO_X(mD)} \ =_{\text{def}} \
\dim_\CC \ \HH{0}{X}{\OO_X(mD)}  \]
as a function of $m$. The precise determination
of these dimensions is of course extremely
subtle even in quite simple situations
\cite{Zariski62b},
\cite{Cutkosky-Srinivas93a},
\cite{Kollar-Matsusaka83}. Our focus will  lie
rather on their asymptotic behavior. In most
interesting  cases the space of   sections in
question grows like
$m^{d}$, and we introduce an invariant that
measures this growth. 
\begin{definition}
\label{Def.of.Volume}
The \textit{volume} of $D$ is
defined to be
\[  \vol_X(D) \ = \  \limsup_{m \to \infty}
\frac{
\hh{0}{X}{\OO_X(mD)}}{m^d/d!}.\]
\end{definition}
\noi The volume of a line bundle is defined
similarly. Note that by definition,
$D$ is big if and only if $\vol_X(D) > 0$. It is
true, but not not entirely trivial, that the
limsup is actually a limit (cf.
\cite[Section 11.4.A]{PAG}). 

In the classical case of \textit{ample}
divisors, the situation is extremely simple. In
fact, if $D$ is ample, then it
follows from the asymptotic Riemann--Roch
theorem
(cf. \cite[1.2.19]{PAG}) that
\begin{align*}
\hh{0}{X}{\OO_X(mD)} \ &= \ \chi\big( X ,
\OO_X(mD) \big) \\ &= \ \big( D^d \big) \cdot
\frac{m^d}{d!} \ + \ O(m^{d-1}) 
\end{align*}
for $m \gg 0$. Therefore, the volume of an
ample divisor $D$ is simply its top
self-intersection number:
\begin{classicaltheorem}[Volume of ample
divisors] \label{Vol.Amp.Div}
If $D$ is ample,
then
\begin{equation} \label{ARR}
\vol_X(D) \ = \ \big ( D^d \big) \ = \ \int_X
c_1\big( \OO_X(D) \big)^d. \notag
\end{equation}
\end{classicaltheorem}
\noi This computation incidentally explains the
terminology: up to constants, the integral in
question computes the volume of $X$ with
respect to a K\"ahler metric arising from the
positive line bundle $\OO_X(D)$. We remark that
the statement remains true assuming only that
$D$ is nef, since in this case
$\hh{i}{X}{\OO_X(mD)} = O(m^{d-1})$ for $i >
0$. 

One can in turn deduce from Theorem
\ref{Vol.Amp.Div} many pleasant features of the
volume function in the classical case. First,
there is a geometric interpretation  springing
from its computation as an intersection
number.
\begin{classicaltheorem}[Geometric
interpretation of
volume] \label{Geom.Interp.Vol.ALB}
 Let $D$ be an ample
divisor. Fix $m
\gg 0$ sufficiently large so that $mD$ is very
ample, and choose $d$ general divisors
\[  E_1, \ldots , E_d \ \in \ \linser{mD}. \]
Then 
\[  \vol_X(D) \ = \ \frac{1}{m^d} \cdot  \#
\, \big( E_1 \cap \ldots \cap E_d \big).\]
\end{classicaltheorem}

Next, one reads off from Theorem
 \ref{Vol.Amp.Div} the behavior of $\vol_X$
as a function of $D$:
\begin{classicaltheorem}[Variational properties
of volume] 
\label{Var.Property.Self.Int}
Given an ample or nef divisor $D$,
$\vol_X(D)$ depends only on the numerical
equivalence class of $D$. It is computed by a
 polynomial function
\[  \vol_X : \textnormal{Nef}(X) \lra \RR   \]
on the nef cone of $X$.
\end{classicaltheorem}

Finally, the intersection form  appearing in
Theorem
\ref{Var.Property.Self.Int} satisfies  a
higher-dimensional extension of the Hodge index
theorem, originally due to Matsusaka, Kovhanski
and Teissier. 
\begin{classicaltheorem}[Log-concavity of
volume]
\label{Log.Concavity.Volume.Classical.Setting}
 Given any two nef classes $\xi, \xi^\pr
\in
\textnormal{Nef}(X)$, one has the inequality
\[\vol_X(\xi + \xi^\pr)^{1/d} \ \ge \
\vol_X(\xi)^{1/d}
\, + \, \vol_X(\xi^\pr)^{1/d}.   
\]
\end{classicaltheorem}
\noi This follows quite easily from 
 the classical Hodge
index theorem on surfaces: see for instance
\cite[Section 1.6]{PAG} and the references cited
therein for the derivation.  
 We refer also to the papers \cite{Gromov90} and
\cite{Okounkov04} of Gromov and Okounkov for an
interesting discussion of this and related
inequalities. 

In the next section, we will see that many of
these properties extend to the case of arbitrary
big divisors. 

\section{Volume of Big Divisors}

We now turn to the volume of arbitrary big
divisors.

It is worth noting right off that there is one
respect in which the general situation differs
from the classical setting.  Namely, it follows
from Theorem
\ref{Vol.Amp.Div} that the volume of an ample
line bundle is a positive integer. However this
is no longer true in general: the volume of a
big divisor can be arbitrarily small, and can
even be irrational.

\begin{example} Let $C$
be a smooth curve, fix
 an integer $a > 1$, and consider the
rank two vector bundle
\[  
E \ = \ \OO_C \big( (1-a) \cdot p \big)\,
\oplus \,
\OO_C\big( p
\big) \]
on $C$, $p \in C$ being some fixed point. Set
\[  X \ = \ \PP(E) \ \ , \ \ L =
\OO_{\PP(E)}(1). \]
Then $L$ is big, and $\vol_X(L) =
\frac{1}{a}$.  (See \cite[Example 2.3.7]{PAG}.)
\end{example}

The first
examples of line bundles having irrational
volume were constructed by Cutkosky
\cite{Cutkosky86a} in the course of his proof of
the non-existence of Zariski decompositions in
dimensions $\ge 3$. We give here a geometric
account of his construction. 
\begin{example}[Integral
divisor with irrational volume]
\label{Int.Div.Irrat.Volume}
  Let $E$ be a general
elliptic curve, and let $V = E \times E$ be the
product of $E$ with itself. Thus $V$ is an
abelian surface, and if $E$ is general then $V$
has Picard number
$\rho(V) =3$, so that $N^1(V)_{\RR} = \RR^3$.
The important fact for us is that
\[ \Effbar(V) \ = \ \Nef(V)
\ \subseteq \ N^1(V) \  = \ \RR^3  \]
is the  circular cone of classes having
non-negative self-intersection (and
non-negative intersection with a hyperplane
class): see Figure \ref{Irrational.Volume}. Now
choose  integral divisors
$A, B$ on
$V$, where
$A$ ample but $B$ is not nef, and write
\[  a \ , b \ \in \  N^1(V)  \]
for the classes of the divisors in question.
Set 
\[   X \ = \ \PP\big(\OO_V(A) \, \oplus \,
\OO_V(B)\big ),
\]
and take $L = \OO_{\PP}(1)$ to be the Serre
line bundle on $X$. We claim that $\vol_X(L)$
will be irrational for general choices of $A$
and $B$. 

\begin{figure}
\hskip .75 in
\includegraphics[scale=.8]{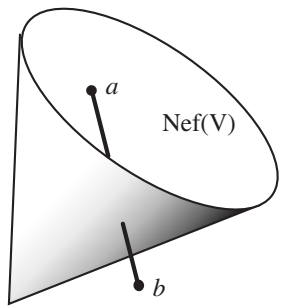}
\caption{Irrational volume
\label{Irrational.Volume}}
\end{figure}

To see this we interpret $\vol_X(L)$ as an
integral. In fact, let   $q(\xi) = \big( \xi
\cdot \xi \big)_V$ be the quadratic intersection
form on $N^1(V)_{\RR}$, and put 
\[  
\widehat{q}(\xi)  \ = \
\begin{cases}
\, q(\xi) & \text{if} \ \xi \in \Nef(V) \\
\,  0  & \text{if} \ \xi \not \in \Nef(V),
\end{cases} 
\]
Then, as in 
\cite[Section 2.3.B]{PAG}, one finds that
\begin{equation}\label{Vol.as.Integral}
\begin{aligned}
\vol_X(L) \ &= \ \frac{3!}{2} \cdot  \int_{0}^1
\widehat{q}\big(  (1-t)a + tb \big) \, dt \\
&= \ \frac{3!}{2} \cdot  \int_{0}^\sigma
 {q}\big(  (1-t)a + tb \big) \, dt,
\end{aligned}
\end{equation}
where $\sigma = \sigma(a,b)$ is the largest
value of $s$ such that $(1-s)a + sb$ is
nef.\footnote{The stated formula  follows via
the isomorphism
\[ \HH{0}{X}{\OO_\PP(m)} \ = \
\HH{0}{V}{ S^m \big( \OO_V(A) \oplus
\OO_V(B) \big) }\] using the fact that     if
$D$ is any divisor on $V$, then 
\[
\hh{0}{V}{\OO_V(D)} \ = \ \begin{cases}
\frac{ (D \cdot D )}{2} & \text{if   $D$
is ample} \\ 0 & \text{if   $D$
is not nef.}
\end{cases}
\]} But for general choices of
$a$ and
$b$,
$\sigma$ is a quadratic irrationality ---
it arises as a root of the quadratic
equation
$q\big( (1-s)a + sb
\big) = 0$ --- so the integral in
\eqref{Vol.as.Integral}    is typically
irrational. The situation is illustrated in
Figure \ref{Irrational.Volume}.
In his thesis \cite{Wolfe05}, Wolfe extends
this interpretation of the volume  as an
integral to much more general projective
bundles. 
\end{example}

We now turn to analogues of the classical
properties of the volume, starting with Theorem
\ref{Geom.Interp.Vol.ALB}. Let $X$ be a smooth
projective variety of dimension $d$, and let
$D$ be a big divisor on $X$. Take a large
integer 
$m
\gg 0$, and fix $d$ general divisors
\[  
E_1, \ldots, E_d \ \in \ \linser{mD}. 
\]  
If $D$ fails to be ample, then it essentially
never happens that the intersection of the
$E_i$ is a finite set. Rather every divisor $E
\in \linser{mD}$ will contain a positive
dimensional \textit{base locus} $B_m \subseteq	
X$, and the $E_i$ will meet along $B_m$
as well as at a finite set of additional points.
One thinks of this finite set as consisting
of the ``moving intersection points" of $E_1,
\ldots , E_d$, since (unlike $B_m$) they vary
with the  divisors $E_i$. 

The following generalization of Theorem
\ref{Geom.Interp.Vol.ALB}, which is essentially
due to Fujita \cite{Fujita94a}, shows that in
general
$\vol_X(D)$ measures the number of these moving
intersection points:

\begin{theorem}[Geometric interpretation of
volume of big divisors] 
\label{Geom.Interp.Volume.Arb.LB}  Always
assuming that $D$ is big,    fix
for $m \gg 0$
$d$ general divisors 
\[ E_1 \,  , \, \ldots \, , \, E_d \ \in \
\linser{mD} \]
when $m \gg 0$, and write $B_m = \Bs{mD}$. Then
\[
\vol_X(D) \ = \ \limsup_{m \to \infty} \frac{\#
\big( E_1
\, \cap \, \ldots \, \cap \, E_d \, \cap \, (X -
B_m)
\big)}{m^d}.
\]
\end{theorem}
\noi The expression in the numerator seems to
have first appeared in work of Matsusaka
\cite{Matsusaka72a},
\cite{Lieberman-Mumford75a}, where it was
called the ``moving self-intersection number"
of $mD$. The statement of
\ref{Geom.Interp.Volume.Arb.LB} appears in
\cite{DEL00}, but it is implicit in Tsuji's
paper \cite{Tsuji99}, and it was
certainly known to Fujita as well. In fact, it
is  a simple consequence of Fujita's theorem,
discussed in \S 3, to the effect that one can
approximate the volume arbitrary closely by the
volume of an ample class on a modification of
$X$.

\begin{remark}[Analytic interpretation of
volume] Fujita's theorem in turn seems to
have been  inspired by some remarks of Demailly
in \cite{Demailly93c}. In that paper, Demailly
decomposed the current corresponding to a
divisor into a singular and an absolutely
continuous part, and he suggested that the
absolutely continuous term should correspond to
the moving part of the linear series in
question. Fujita seems to have been lead to his
statement by the project of algebraizing
Demailly's results  (see also \cite{ELN96} and
\cite[\S 7]{Lazarsfeld97a}). The techniques and
intuition introduced by Demailly have been
carried forward  by Boucksom and others
\cite{Boucksom02}, \cite{Boucksom04},
\cite{BDPP} in the analytic approach to
asymptotic invariants. 
\end{remark}

We turn next to the variational properties of
$\vol_X$. Observe to begin with that $\vol_X(D)$
is naturally defined for any $\QQ$-divisor $D$.
Indeed, one can  adapt Definition
\ref{Def.of.Volume} by taking the limsup over
values of $m$ that are sufficiently divisible
to clear the denominators of $D$.
Alternatively, one can check (cf.
\cite[2.2.35]{PAG}) that the volume on integral
divisors satisfies the homogeneity property
\begin{equation} \label{Homogeneity.Volume}
\vol_X(aD) \ = \ a^d \cdot \vol_X(D),
\end{equation}
and one can use this in turn to define
$\vol_X(D)$ for $D \in \Div_{\QQ}(X)$.

The following results were proved by the second
author in \cite[Section
2.2.C]{PAG}. They are the analogues for
arbitrary big divisors of   the classical
 Theorem \ref{Var.Property.Self.Int}.
\begin{theorem}[Variational properties of
$\vol_X$] 
\label{Variational.Properties.Volume}
Let $X$ be a smooth
projective variety of dimension $d$.
\begin{enumerate} 
\item[(i)] The volume of a   $\QQ$-divisor on
$X$ depends only on its numerical equivalence
class, and hence this invariant defines a
function
\begin{equation}  \vol_X  \, : \, N^1(X)_{\QQ}
\lra
\RR. \notag \end{equation}
\item[(ii).] 
Fix any norm
$\Vert \  \Vert$ on
$N^1(X)_{\RR}$ inducing the usual topology on 
that   vector space. Then
there is a positive constant $C > 0$ such that
\begin{equation} \label{contin.vol.eqn}
\big \vert \, \vol_X\big( \xi \big) \, -
\,\vol_X
\big(
\xi^\pr \big) \, \big \vert \ \le
\ C \cdot \Big( \max \big( \, \Vert \xi\Vert\, 
, \, 
\Vert \xi^\pr \Vert\, \big) \Big)^{d -1} \cdot
\Vert \,
\xi \, - \, \xi^\pr \, \Vert  \notag
\end{equation}
for any two classes $\xi, \
\xi^\pr \in N^1(X)_{\QQ}$. 
\end{enumerate}
\end{theorem}

\begin{corollary}[Volume of real
classes]
\label{vol.real.classes}
  The function $\xi \mapsto
\vol_X(\xi)$ on
$N^1(X)_\QQ $ extends uniquely to a continuous function 
\[ \vol_X \, : \,  N^1(X)_{\RR} \lra \RR. \qed
\]
\end{corollary}
\noi 
We note that both the theorem and its corollary
hold for arbitrary irreducible projective
varieties.  For smooth complex manifolds, the
continuity of volume was established
independently by Boucksom in \cite{Boucksom02}.
In fact, Boucksom defines and studies the volume
of an arbitrary pseudoeffective class $\alpha
\in
\HH{1,1}{X}{\RR}$ on a compact  K\"ahler
manifold $X$: this involves some quite
sophisticated analytic methods. 

One may say that the entire asymptotic 
Riemann--Roch problem on $X$ is encoded in the
finite-dimensional vector space
$N^1(X)_{\RR}$ and the continuous function
$\vol_X : N^1(X)_{\RR} \lra \RR$ defined on
it. So it seems like a rather basic problem to  
understand the behavior of this invariant as
closely as possible.

\begin{figure}
\includegraphics[scale=1]{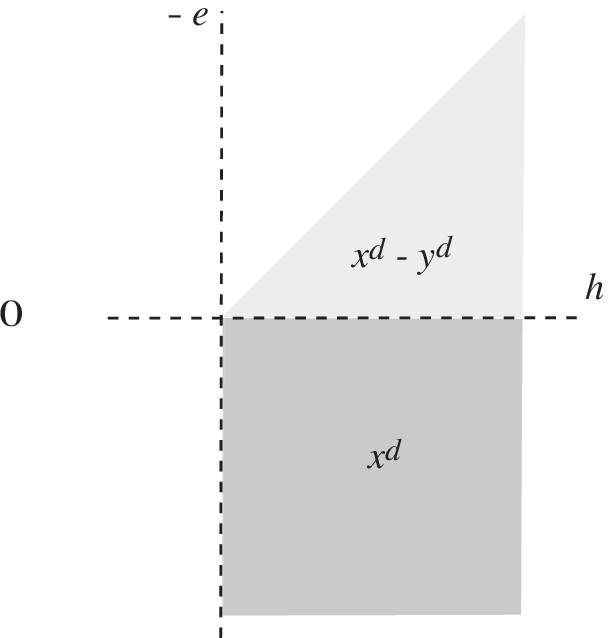}
\caption{Volume on blow-up of $\PP^d$
\label{Volume.On.BlowUp}}
\end{figure}

We next present some examples and
computations.
\begin{example} 
[Volume on blow-up of projective space]
\label{Volume.BlowUp.Proj.Space.Ex}
The first non-trivial example
is to take $X = \Bl_p{\PP^d}$ to be the
blowing-up of projective space $\PP^d$ at a
point $p$. Write $H$ for the pull-back of a
hyperplane, and $E$ for the exceptional
divisor, and denote by $h, e \in N^1(X)_{\RR}$
  the corresponding classes. Then
\[  N^1(X)_{\RR} \ = \ \RR \cdot h + \RR \cdot
e. \]
 The nef cone $\Nef(X)$ is
generated by $h-e$ and $h$, while the
pseudoeffective cone is spanned by $h-e$ and
$e$. On $\Nef(X)$, the volume is just given by
the intersection form:
\[   \vol_X( xh - ye) \ = \ \big(  (xh - ye)^ d
\big) \ =
\ x^d - y^d \ \ \ (0 \le x \le y).
\]
The sector of $N^1(X)$ spanned by $h$ and
$e$ corresponds to linear series of the form
$\linser{aH + bE}$ with $b \ge 0$, and one
checks that the exceptional divisor is fixed in
those linear series, i.e.
\[ \HH{0}{X}{\OO_X(aH + bE)} \ = \ 
\HH{0}{X}{\OO_X(aH)}\]
when $a, b \ge 0$. Thus if $x \ge 0 \ge y$,
then \[
\vol_X(xh -ye)  \ = \ \vol_X(xh) \ = \ x^d.\]
Elsewhere the  $\vol_X = 0$. The situation is
summarized in Figure \ref{Volume.On.BlowUp},
which shows
$\vol_X(xh -ye)$ as a function of $(x,y) \in
\RR^2$. 
 \end{example}

\begin{example}[Toric varieties]  When $X$ is a
toric variety, there is a
fan refining the pseudoeffective cone of $X$,
the
\emph{Gel'fand-Kapranov-Zelevinski
decomposition} of \cite{Oda-Park91}, with
respect to which the volume function is
piecewise polynomial. In fact, recall that every
torus-invariant divisor on $X$ determines a
polytope $P_D$ such that
$\hh{0}{X}{\OO_X(D)}$ is the number of
lattice points in $P_D$ and $\vol_X(D)$ is
the lattice volume of $P_D$. On the interior of
each of the maximal cones in the GKZ
decomposition the combinatorial type of the
polytope $P_D$ is constant, i.e.   all the
polytopes have the same normal fan $\Sigma_D$,
$D$ is the pull-back of an ample divisor $A$ on
the toric variety $Y$ corresponding to
$\Sigma_D$, and ${\rm vol}_X(D)={\rm vol}_Y(A)$.

The family of polytopes $\{P_D\}_D$ can be
considered as a
\emph{family of partition polytopes} and the
function
$D\mapsto \hh{0}{X}{\OO_X(D)}$ is the
corresponding
\emph{vector partition function}. Brion and
Vergne \cite{Brion-Vergne97}   have studied  
the continuous function associated to a
partition function (in our case this is
precisely the volume function), and they gave
explicit formulas on each of the chambers in
the corresponding fan decomposition.
\end{example}

As we have just noted, if $X$ is a toric
variety then
$\vol_X$ is piecewise polynomial with respect
to a polyhedral subdivision of $N^1(X)_{\RR}$.
This holds more generally when   the
linear series on
$X$ satisfy a very strong hypothesis of finite
generation.
\begin{definition} [Finitely generated linear
series]
\label{Finite.Gen.Lin.Series.Def}
 We say that $X$ has \textit{finitely
generated linear series} if there exist
integral divisors
$D_1, \ldots , D_r$ on $X$, whose classes form
a basis of $N^1(X)_{\RR}$, with the property
that the $\ZZ^r$-graded ring
\[ R\big(X; D_1, \ldots , D_r) \ = \
\bigoplus_{m_1, \ldots m_r \in \ZZ}
\HH{0}{X}{\OO_X(m_1D_1 + \ldots + m_rD_r)}
\] is finitely generated.
\end{definition}
\noi 
One should keep in mind that  the finite
generation of this ring does not in general
depend only on the numerical equivalence
classes of the $D_i$. The definition was
inspired by a very closely related concept
introduced and studied by Hu and Keel in
\cite{Hu-Keel00}.

\begin{example} It is a theorem of Cox
\cite{Cox95} that the condition in the
definition is satisfied when
$X$ is a smooth toric variety, and
\ref{Finite.Gen.Lin.Series.Def} also holds when
$X$ is a non-singular  \textit{spherical
variety} under a reductive group.\footnote{As
M. Brion kindly explained to us, the finite
generation for spherical varieties follows
from a theorem of Knop in \cite{Knop93}.}  It is
conjectured in
\cite{{Hu-Keel00}} that smooth Fano varieties
also have finitely generated linear series:
this is verified by Hu and Keel in dimension
three using the minimal model program. 
\end{example}

Again inspired by the work of Hu and Keel just
cited, it was established by the authors in
\cite{AIBL} that the piecewise polynomial nature
of the volume function observed in Example
\ref{Volume.BlowUp.Proj.Space.Ex} holds in
general on varieties with finitely generated
linear series.
\begin{theorem}
Assume that $X$ has finitely generated linear
series. Then $N^1(X)_{\RR}$ admits a finite
polyhedral decomposition with respect to which
$\vol_X$ is piecewise polynomial.
\end{theorem}
\noi Altough the method of proof in \cite{AIBL}
is different, in the most important cases one
could deduce the theorem directly from the
results of
\cite{{Hu-Keel00}}.

\begin{figure}
\vskip 10pt
\includegraphics[scale=.9]{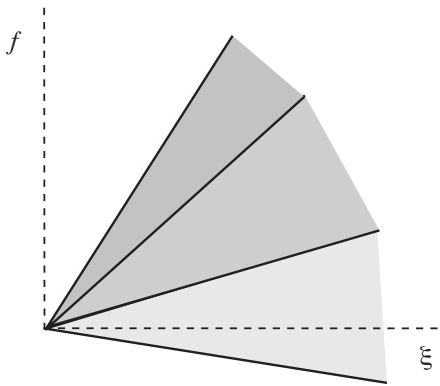}
\caption{Volume on projective bundle over curve
\label{Volume.On.Bundle.Curve}}
\end{figure}

\begin{example} [Ruled varieties over curves] 
Wolfe  \cite{Wolfe05} analyzed the volume
function on ruled varieties over curves, and
found a pleasant connection with some classical
geometry. Let
$C$ be a smooth projective curve of genus $g
\ge 1$, and $E$ a vector bundle on $C$ of rank
$e$. We consider  $X = \PP(E)$. Then
\[  N^1(X)_{\RR} \ = \ \RR \cdot  \xi \, + \,
\RR \cdot f, \]
where $\xi = c_1\big(\OO_{\PP(E)}(1)\big)$ is
the class of the Serre line bundle and $f$ is
the class of a fibre of the bundle map $\PP(E)
\lra C$. There are now two cases:
\begin{itemize}
\item If $E$ is
\textit{semistable},\footnote{This means roughly
speaking that $E$ does not admit any quotients
of too small degree. See for
instance \cite{Huybrechts-Lehn97} or
\cite[Section 6.4.A]{PAG}.}  then it is an old
observation of Miyaoka that 
\[\Nef(X) \ = \ \Effbar(X). \]
So in this case $\vol_X$ is just given by the
intersection form on the nef cone.
\lbl
\item When $E$ is not semistable, it admits a
canonical \textit{Harder-Narasimhan filtration} 
with semistable graded pieces. Wolfe shows
that this determines a decomposition of
$\Effbar(X)$ into finitely many sectors --- one
for each piece of the filtration --- on which
$\vol_X$ is given by a
polynomial.\footnote{Interestingly, it has not
so far proved practical to  evaluate
these polynomials explicitly.} 
\end{itemize}
The situation is illustrated in Figure
\ref{Volume.On.Bundle.Curve}. In summary we may
say that the basic geometry associated to the
vector bundle $E$ is visible in the volume
function on $X = \PP(E)$. 
\end{example}

\begin{problem}
It would be worthwhile to work out examples
of the volume function $\vol_X$ when $X$ is
the projectivization of a bundle $E$ on a higher
dimensional base. Here one might expect more
complicated behavior, and it would be
interesting to see to what extent the geometry
of $E$ is reflected in $\vol_X$.  
\end{problem}

In the case of surfaces, the volume function
was analyzed by Bauer, K\"uronya and Szemberg
in \cite{BKS}, who prove:
\begin{theorem}[Volume function on surfaces] Let
$X$ be a smooth surface. Then
$\textnormal{Big}(X)$ admits a locally finite
rational polyhedral subdivision on which
$\vol_X$ is piecewise polynomial.
\end{theorem}
\noi The essential point here is that on
surfaces, the volume of a big divisor $D$ is
governed by the \textit{Zariski decomposition}
of $D$, which gives a canonical expression $D =
P + N$ of the divisor in question as the sum of
a ``positive" and ``negative part" (see for
instance
\cite{Badescu01a}, or \cite{BKS}).
It is immediate that $\vol_X(D) = \big(
P^2\big)_X$, so the issue is to understand how
this decomposition varies with $D$. 

On the other hand, in the same paper
\cite{BKS},  Bauer, K\"uronya and Szemberg give
examples to prove:
\begin{proposition} There exists a smooth
projective threefold $X$ 
  for which \[
\vol_X : N^1(X)_{\RR} \lra \RR\] is    
piecewise analytic but not piecewise polynomial.
\end{proposition}
\noi This shows that in general, the volume
function has a fundamentally non-classical
nature.

\begin{example}[Threefold with non-polynomial
volume] The idea of \cite{BKS} is to study the
volume function on the threefold $X$
constructed in Example
\ref{Volume.BlowUp.Proj.Space.Ex}. Keeping the
notation of that example, fix a divisor $C$ on
$V$ with class $c \in N^1(X)_{\RR}$. In view of
the homogeneity \eqref{Homogeneity.Volume} of
volume, the essential point is to understand
how the volume of the line bundle
\[L_C \ =_{\text{def}} \  \OO_{\PP}(1)
\otimes
\pi^*
\OO_V(C)
\] varies with $C$, where $\pi : X
= \PP\big(\OO_V(A)
\oplus \OO_V(B) \big)  \lra V$ is the bundle
map.\footnote{In the first instance $C$ is an
integral divisor, but the computation that
follows works as well when $C$ has $\QQ$ or
$\RR$ coefficients.} For this one can proceed as
in Example
\ref{Volume.BlowUp.Proj.Space.Ex}.
One finds that 
\[  \vol_X(L_C) \ = \ \frac{3!}{2}  \cdot 
\int_{0}^{\sigma(c)}
 {q}\big( c + (1-t)a + tb \big) \, dt,\]
where $\sigma(c)$ denotes the largest value of
$s$ for which
$(1-s)a + sb +c$ is
nef.\footnote{Geometrically,
$\sigma(c)$ specifies the point at which the
line segment joining $a +c$ to $b + c$ meets the
boundary of the circular nef cone of $V$.}
On the other hand, if $c$ is small
enough so that $a +c$ is nef while $b+c$ is
not, then $\sigma(c)$ is not given by a
polynomial in 
$c$ (although
$\sigma(c)$ is an algebraic function of $c$).
So one does not expect --- and it is not the
case --- that
$\vol_X(L_C)$ varies polynomially with $C$.  
\end{example}

Returning to the case of a smooth projective
variety $X$ of arbitrary dimension $d$, the
second author observed  in \cite[Theorem
11.4.9]{PAG} that  Theorem
\ref{Log.Concavity.Volume.Classical.Setting}
extends to arbitrary big divisors:
\begin{theorem}[Log-concavity of volume] 
\label{Log.Concavity.Volume}
The
inequality
\[\vol_X(\xi + \xi^\pr)^{1/d} \ \ge \
\vol_X(\xi)^{1/d}
\, + \, \vol_X(\xi^\pr)^{1/d}.   
\]
holds for any two classes $\xi, \xi^\pr \in
\textnormal{Big}(X)$. 
\end{theorem}
\noi Although seemingly somewhat delicate, this
actually follows immediately from the theorem
of Fujita that we discuss in the next section.

As of this writing, Theorems
\ref{Variational.Properties.Volume} and
\ref{Log.Concavity.Volume} represent the only
regularity properties that the volume function
$\vol_X$ is known to satisfy in general. It
would be very interesting to know whether in
fact there are others. The natural expectation
is that
$\vol_X$ is ``typically" real analytic:
\begin{conjecture}
There is a ``large" open set $U \subseteq
\textnormal{Big}(X)$ such that $\vol_X$ is real
analytic on each connected component of $U$.
\end{conjecture}
\noi One could hope more precisely that $U$ is
actually dense, but this might run into trouble
with clustering phenomena. We refer to Section
6 for a description of some ``laboratory
experiments" that have been performed on a
related invariant.

\section{Fujita's Approximation Theorem}

One of the most important facts about $\vol_X$
--- and the source of several of the results
stated in the previous section --- is a
theorem of Fujita \cite{Fujita94a} to the effect
that one can approximate the volume of an
arbitrary big divisor by the top-self
intersection number of an ample line bundle on
a blow-up of the variety in question. We
briefly discuss Fujita's result here.

As before, let $X$ be a smooth projective
variety and consider a big class $\xi \in
N^1(X)_{\RR}$. As a matter of terminology, we
declare that a 
\textit{Fujita approximation} of
$\xi$ consists of a birational morphism
\[  \mu : X^\pr \lra X, \] from a smooth
projective variety $X^\pr$ onto $X$, together
with a decomposition
\[  \mu^* \xi \ = \ a \, + \, e \ \in \
N^1(X^\pr)_{\RR},\]
where $a$ is an ample class  and $e$
is an effective class on $X^\pr$.\footnote{By
definition, an effective class in
$N^1(X^\pr)_{\RR}$ is one which is represented
by an effective
$\RR$-divisor.}

Fujita's theorem asserts that one can find such
an approximation in which the volume of $a$
approximates arbitrarily closely the volume of
$\xi$:
\begin{theorem}[Fujita's approximation theorem]
\label{Fujita.Approx.Thm}
Given any $\eps > 0$, there exists a Fujita
approximation 
\[  \mu : X^\pr \lra X \ \ , \ \ \mu^*(\xi) \ =
\ a \, + \, e \]
where $ \vert \vol_X(\xi) - \vol_{X^\pr}(a)
\vert < \eps$.
\end{theorem}
\noi We stress that the output of the theorem
depends on $\eps$. We note that Theorems
\ref{Geom.Interp.Volume.Arb.LB}  and
\ref{Log.Concavity.Volume} follow almost
immediately  from this result.

For the proof, one reduces by continuity and
homogeneity to the situation when $\xi$
represents an integer divisor $D$. Given $m \gg
0$, denote by $\frb_m \subseteq \OO_X$ the
base-ideal of $\linser{mD}$, and let
$
\mu_m : X_m \lra X 
$ be a resolution of singularities of the
blowing-up of $\frb_m$.  Then on $X_m$
one can write
\[  \mu_m^*(mD) \ = \ B_m \, + \, E_m ,
\]
where $B_m$ is free and $E_m$ is (the  pullback
of) the  exceptional divisor of the blow-up. 
The strategy is to take $A_m = \tfrac{1}{m}
B_m$ as the positive part of the
approximation.\footnote{$A_m$ is not
 ample, but it is close to being
so.} It is automatic that $\vol_{X_m}(A_m) \le
\vol_X(D)$, so the essential point is to prove
that
\[ \vol_{X_m}(A_m) \ > \ \vol_X(D) \, - \, \eps 
\] provided that $m$ is sufficiently large
(with the bound depending on $\eps$). Fujita's
original proof was elementary but quite tricky;
a more conceptual proof using multiplier ideals
appears in
\cite{DEL00}, and yet another argument, 
revolving around higher jets, is given in
\cite{Nakamaye03}. 

Boucksom, Demailly, Paun and Peternell
established in \cite{BDPP} an important
``orthogonality" property of Fujita
approximations. Specifically, they prove:
\begin{theorem}
In the situation of Theorem
\ref{Fujita.Approx.Thm}, fix an ample class $h$
on $X$ which is sufficiently positive so that
$h \pm \xi$ is ample. Then there is a
universal constant $C$ such that
\[\big( a^{n-1} \cdot e \big)^2_{X^\pr} \ \le \
C \cdot\big(h^n\big)_X \cdot \big( \vol_X(\xi)
\, -
\,
\vol_{X^\pr}(a) \big). \]
\end{theorem}
 \noi In
other words,
as one passes to a better and better Fujita
approximation, the exceptional divisor of the
approximation becomes more nearly orthogonal to
the ample part of $\xi$. This had been
conjectured by the fourth author in
\cite{Nakamaye03}. 
Boucksom--Demailly--Paun--Peternell use this to
identify the cone of curves dual to the
pseudoeffective cone. Accounts of this work
appears in \cite{Debarre05} and in \cite[\S
11.4.C]{PAG}.

\section{Higher Cohomology}

The volume $\vol_X(D)$ measures the asymptotic
behavior of the
dimension $h^0\big(mD\big)$. But of course
it has been well-understood since Serre that
one should also consider higher cohomology
groups. We discuss in this section invariants
involving the groups $\HH{i}{X}{\OO_X(mD)}$ for
$i \ge 1$. 

As always, $X$ denotes a smooth projective
variety of dimension $d$, and $D$ is a divisor
on $X$.
\begin{definition} [Asymptotic cohomology
function] Given $i \ge 0$, the $i^{\text{th}}$
\textit{asymptotic cohomology function}
associated to
$D$ is
\[ \hhat{i}{D} \ = \
\limsup
\,
\frac{\hh{i}{X}{\OO_X(mD)}}{m^d/d!}.\]
\end{definition}
\noi Thus $\vol_X(D) = \hhat{0}{D}$.  By working with sufficiently
divisible $m$, or establishing the homogeneity
of $\widehat{h}^i$ as in equation
\eqref{Homogeneity.Volume}, the definition
extends in a natural way to
$\QQ$-divisors. However we remark that it is
not known that the limsup in the definition is
actually a limit, although one hopes that this
is the case.

These higher cohomology functions were studied 
in \cite{Kuronya05}   by K\"uronya,
who establishes for them the analogue of Theorem
\ref{Variational.Properties.Volume}:
\begin{theorem}
The function $\hhat{i}{D}$ depends only on the
numerical equivalence class of a $\QQ$-divisor
$D$, and it extends uniquely to a continuous
function 
\[  \widehat{h}^i : N^1(X)_{\RR} \lra \RR. \]
\end{theorem}
\noi Note that in contrast to the volume,
these functions can be non-vanishing away from
the big cone. For example, it follows from 
Serre duality that the highest function
$\widehat{h}^d$ is supported on $-
\Effbar(X)$.

\begin{example} [Blow-up of $\PP^d$] Returning
to the situation of Example
\ref{Volume.BlowUp.Proj.Space.Ex}, let $X =
\Bl_p(\PP^d)$ be the blowing up of projective
space at a point, and consider the region $x >
0 > y$
of $N^1(X)_{\RR} = \RR^2$ where the class $xh -
ye$ is big but not ample. In this sector,
$\widehat{h}^i = 0$ unless $i =0 $ or $i = d-1$,
and 
\[  \widehat{h}^0\big(xh
- ye\big) = x^d \ \ \ ,
\ \ \ \ \widehat{h}^{d-1}\big(xh - ye
\big)
\ =
\ (-1)^d \cdot y^d. \]

\end{example}

\begin{figure}
\vskip 10pt
\includegraphics[scale=.8]{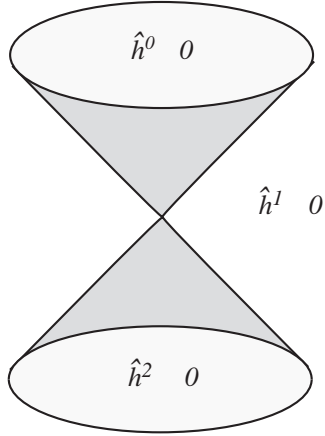}
\caption{Higher cohomology on abelian surface
\label{Higher.Cohom.Abel.Sf.Figure}}
\end{figure}

\begin{example} [Abelian varieties,
\cite{Kuronya05},
\S 3.1] 
\label{Higher.Cohom.Abel.Vars.Example}
Let
$X$ be an abelian variety of dimension $d$.
Here the  higher cohomology functions reflect
the classical theory of indices of line bundles
on abelian varieties.  Specifically,  express
$X$ as the quotient
$X = V/
\Lambda$ of a
$d$-dimensional complex vector space $V$  by a
lattice $\Lambda \subseteq V$. Then
$\HH{1,1}{X}{\RR}$ is identified with the real
vector space of hermitian forms on $V$, and
$N^1(X)_{\RR} \subseteq \HH{1,1}{X}{\RR}$ is the
subspace spanned by those forms whose imaginary
part takes integer values on
$\Lambda$. The dense open subset of this
N\'eron--Severi space corresponding to
non-degenerate forms is partitioned into
disjoint open cones  
\[
\textnormal{Ind}_j(X) \ \subseteq \
N^1(X)_{\RR}
\]
 according to the
\text{index} (i.e. the number of negative
eigenvalues) of the form in question. Then given
$\xi \in N^1(X)_{\RR}$, one has
\[
\hhat{i}{\xi} \ = \ 
\begin{cases}
(-1)^i \cdot \int_X \xi^d & \text{ if } \, \xi
\,
\in\,
\textnormal{Ind}_i(X) \\
0 & \text{ if } \, \xi\, \not \in \,
\textnormal{Ind}_i(X).
\end{cases}
\]
For
instance, suppose  
that $X$ has dimension $2$ and Picard
number $\rho(X) = 3$. Then (as in Example
\ref{Int.Div.Irrat.Volume}) the cone of classes
having self-intersection $ = 0$ divides
$N^1(X)_{\RR}$ into three regions. The classes
with $\widehat{h}^0 \ne 0$ and $\widehat{h}^2
\ne 0$ occupy the two conical regions, while
the exterior of the cone consists of classes
with $\widehat{h}^1 \ne 0$. The situation is
illustrated in Figure
\ref{Higher.Cohom.Abel.Sf.Figure}.
\end{example}

\vskip 10pt

\begin{example}[Flag varieties] Let
$X = G/B$ be the quotient of a semi-simple
algebraic group by a Borel subgroup. Then the
 Borel--Bott--Weil theorem gives a finite
 decomposition of
$N^1(X)_{\RR}$ into polyhedral chambers on which
exactly one of the functions $\widehat{h}^i$ is
non-zero. See \cite[\S 3.3]{Kuronya05} for
details. 
\end{example}

\begin{remark}
[Toric varieties] When $X$ is a projective
toric variety, the functions $\widehat{h}^i$
are studied by Hering, K\"uronya and Payne in
\cite{HKP05}. In this case, the functions are
again piecewise polynomial on a polyhedral
decomposition of $N^1(X)_{\RR}$, although in
general one needs to refine the
GKZ-decomposition with respect to which $\vol_X$
is piecewise polynomial. 
\end{remark}

It is interesting to ask what geometric
information these functions convey.   De
Fernex, K\"uronya and the second author
\cite{DKL05} have established that ample classes
are characterized by vanishings involving the
higher cohomology functions.
\begin{theorem} [Characterization of ample
classes] 
\label{Theorem.Of.DKL}
Let \[ \xi_0 \ \in\ 
N^1(X)_{\RR}\] be any class. Then $\xi_0$ is
ample if and only if there exists a
neighborhood $ U
\subseteq N^1(X)_{\RR}$ of $\xi_0$ such that
\[  \widehat{h}^i (\xi) \ = \ 0 \]
for all $i > 0$ and every $\xi \in U$.
\end{theorem}
 \noi As one sees for instance from Example
\ref{Higher.Cohom.Abel.Vars.Example}, it can
happen that 
$\widehat{h}^i(\xi_0) = 0$  for
all $i$ without $\xi_0$ being ample (or even
pseudoeffective), so one can't avoid looking at
a neighborhood of
$\xi_0$. One can view Theorem
\ref{Theorem.Of.DKL} as an asymptotic analogue
of Serre's criterion for amplitude. The
essential content of the result is that if $D$
is a big divisor that is not nef, and if $A$
is any ample divisor, then $\hhat{i}{D -tA}
\ne 0$ for some $i > 0$ and sufficiently small
$t > 0$.

\section{Base Loci}

Big divisors that fail to be nef are
essentially characterized by the asymptotic
presence of base-loci. It is then natural to
try to measure quantitatively the loci in
question. This was first undertaken by Nakayama
\cite{Nakayama04}, and developed from another
viewpoint in the paper \cite{AIBL} of the
present authors. Nakayama's results were
generalized and extended to the analytic
setting by Boucksom in \cite{Boucksom04}.

As before, let $X$ be a smooth projective
variety of dimension $d$. Denote by $ 
\ord_E$ a divisorial valuation centered on $X$.
Concretely, this is given by a projective
birational morphism $\mu : X^\pr \lra X$, say
with $X^\pr$ smooth, together with an
irreducible divisor $E \subseteq
X^\pr$.\footnote{Note that the same valuation
can arise from different models
$X^\pr$.} Given locally a regular function $f$
on $X$, one can then discuss the order
\[  \ord_E(f) \ \in \ \NN \]
of vanishing of  $f$ along $E$. For an
effective divisor $D$ on $X$, the order
$\ord_E(D)$ of $D$ along $E$ is defined via a
local equation for $D$. 

\begin{example} [Order of vanishing at a
point or along a subvariety] Fix a point
$x
\in X$, and let
$E$ denote the exceptional divisor in the
blow-up $\Bl_x(X)$. Then $\ord_E = \ord_x$ is
just the classical order of vanishing at $x$:
one has 
\[ \ord_x(f) \ = \ a \]
if and only if all partials of $f$ of order $<
a$ vanish at $x$, while some $a^{\text{th}}$
partial is $\ne 0$. Given any proper 
irreducible subvariety $Z \subseteq X$, the
order of vanishing $\ord_Z$ is defined
similarly. 
\end{example}

Now let $D$ be a big divisor on $X$. Then one
sets
\begin{align*} 
\ord_E\big( \linser{D} \big) \ &= \ \min_{D^\pr
\in \linser{D}} \ord_E \big( D^\pr \big) \\
\ord_E\big( \alinser{D} \big) \ &= \ \limsup
\frac{
\ord_E \big(\linser{mD} \big)  }{m}.
\end{align*}
These definitions once again extend in a
natural way to $\QQ$-divisors, and then one has:
\begin{theorem}
The function $\ord_E\big( \alinser{D} \big)$
depends only on the numerical equivalence class
of a big $\QQ$-divisor $D$, and it extends to a
continuous function
\begin{equation}
\label{Order.Function}  \ord_E \,  : \,
\textnormal{Big}(X) \lra \RR.
\end{equation}
\end{theorem}
\noi This was initially established by Nakayama
for order of vanishing along a subvariety, and
extended to more general valuations in
\cite{AIBL}.

The paper \cite{AIBL} also analyzes precisely 
when this invariant vanishes. Recall for this
that the \textit{stable base-locus} $\SBL{D}$
of a big divisor (or $\QQ$-divisor) $D$ is
defined to be the common base-locus of the
linear series $\linser{mD}$ for sufficiently
large and divisible $m$. (See
\cite[Section 2.1.A]{PAG} for details.)
\begin{theorem} Keeping notation as above, let
 $Z$ be the center of $E$ on $X$. If $D$ is a
big  divisor on $X$, then 
\[ \ord_E\big( \alinser{D} \big) \ > \ 0\]
if and only if there exists an ample
$\QQ$-divisor $A$ such that $Z$ is  
contained in the stable base-locus $\SBL{D +
A}$ of $D + A$.  
\end{theorem}
\noi One thinks of $D +A$ as a small
positive perturbation of $D$. The union of
the sets
$\SBL{D + A}$ is called in \cite{BDPP} the
\textit{non-nef} locus of $D$. With a little
more care about the definitions, one can prove
the analogous statement with $D$ replaced by an
arbitrary big class $\xi \in
\text{Big}(X)$. 

Building on ideas first introduced by the
fourth author in \cite{Nakamaye01a} and
\cite{Nakamaye03}, one can also use a
volume-like invariant to analyze what
subvarieties appear as irreducible components
of  the stable base-locus. In fact, let
$Z
\subseteq X$ be an irreducible subvariety of
dimension $e$, and let $D$ be a big divisor on
$X$. Set
\[  V_m \ = \ \text{Image} \Big (
\HH{0}{X}{\OO_X(mD)}  \lra
\HH{0}{Z}{\OO_Z(mD)} \Big). \]
 Thus $V_m$ is a (possibly proper) subspace of
the space of sections of $\OO_Z(mD)$.
We define the
\textit{restricted volume} of $D$ to $Z$ to be
\[ 
\vol_{X | Z}(D) \ = \ \limsup \frac{ \dim
V_m}{m^e / e!}. 
\]
For instance if $D$ is ample then the
restriction mappings are eventually surjective,
and hence \[
\vol_{X | Z}(D) \ = \ \vol_Z(D\vert Z) \ = \ 
\big( D^e
\cdot Z
\big). \] However in general it can happen
that $\vol_{X | Z}(D) < \vol_Z(D\vert Z)$. An
analogue of Theorem
\ref{Geom.Interp.Volume.Arb.LB} holds for these
restricted volumes. Again the definition
extends naturally to $\QQ$-divisors
$D$ (and $\RR$-divisors as well).  

Given any $\QQ$-divisor $D$, the
\textit{augmented base-locus} of $D$ is defined
to be
\[  \BAug{D} \ = \ \SBL{D -   A} \]
for any ample $\QQ$-divisor $A$ whose class in
$N^1(X)_{\RR}$ is sufficiently small (this being
independent of the particular choice of
$A$).  It was discovered by the
fourth author in
\cite{Nakamaye01a} that this augmented
base-locus behaves more predictably that
$\SBL{D}$ itself. (See
\cite[Chapter 10.4]{PAG} for an account.)

Building on the ideas and techniques of
\cite{Nakamaye03}, the authors prove in
\cite{ELMNP2}
\begin{theorem}
Let $D$ be a  big $\QQ$-divisor on $X$, and let
$Z \subseteq X$ be a subvariety. If $Z$
is an irreducible component of
$\BAug{D}$, then
\[  \lim_{D^\pr \lra D} \vol_{X | Z}(D^\pr) \
=
\ 0, \]
the limit being taken over all big
$\QQ$-divisors approaching $D$.
Conversely, if $Z \not \subseteq \BAug{D}$,
then 
\[  \lim_{D^\pr \lra D} \vol_{X \vert Z}(D^\pr) \
= \ \vol_{X \vert Z}(D) \ > 
\ 0. \]
\end{theorem}
\noi We remark that  the statement
remains true also when $D$ is an $\RR$-divisor.

\section{In Vitro Linear Series}

Just as it is interesting to ask about the
regularity of the volume function, it is also
natural to inquire about the nature of the
order function occuring in
\eqref{Order.Function}. To shed light on this
question, one can consider an abstract algebraic
construction that models some of the global
behavior that can occur for global linear
series.

Let $X$ be a smooth projective variety of
dimension $d$. The starting point is the
observation that the pseudoeffective and nef
cones 
\[   \Nef(X) \ \subseteq \
\Effbar(X)
\ \subseteq \ N^1(X)_{\RR}
  \]
on $X$, 
as well as the order function in
\eqref{Order.Function}, can be recovered
formally from the base-ideals of   linear
series on $X$. Specifically, fix integer
divisors $D_1, \ldots , D_r$ on $X$ whose
classes form a basis of $N^1(X)_{\RR}$: note
that this determines an identification
$N^1(X)_{\RR} = \RR^r$. Given inegers $m_1,
\ldots, m_r \in \ZZ$, write $\vec{m} = (m_1,
\ldots , m_r)$, and let
\[ \frb_{\vec{m}} \ = \
\bs{ m_1D_1 + \ldots + m_rD_r}  \] 
be the base ideal of the corresponding linear
series. Thus $\frb_{\vec{m}} \subseteq \OO_X$
is an ideal sheaf on $X$, and one has
\begin{equation}  \label{Mult.Base.Ideals}
\frb_{\vec{m}} \cdot
\frb_{\vec{\ell}}
\ \subseteq \ \frb_{\vec{m} +
\vec{\ell}}\end{equation} for every $ \vec{m},
\vec{\ell} \in \ZZ^r$.  Note that we can
recover the nef and pseudoeffective cones on
$X$ from these ideals. In fact, a moment's
thought shows that
$\Nef(X) \subseteq \RR^r$ is the closed convex
cone spanned by all vectors $\vec{m} \in
\ZZ^r$ such that $\frb_{\vec{m}} = \OO_X$, and
similarly $\Effbar(X)$ is the closed
convex cone generated by all
$\vec{m} $ such that $\frb_{\vec{m}} \neq
(0)$. The  function 
\[ \ord_Z \,  : \,
\textnormal{Big}(X)_{\RR}\lra \RR
\] measuring order of vanishing along a
subvariety $Z \subseteq X$ discussed in the
previous section can also be defined just using
these ideals. 

This leads one to consider arbitrary
collections of ideal sheaves satisfying the
basic multiplicativity property
\eqref{Mult.Base.Ideals}: they provide abstract
models of linear series. Specifically, let $V$
be any smooth variety of dimension $d$ (for
example $V = \CC^d)$. A \textit{multi-graded 
family of ideals} on $V$ is a collection 
\[  \fra\bull \ = \ \big \{ \fra_{\vec{m}} \big
\}_{\vec{m} \in \ZZ^m} \]   of
ideal sheaves $\fra_{\vec{m}} \subseteq \OO_V$
on
$V$, indexed by
$\ZZ^r$, such that  $\fra_{\vec{0}} = \OO_V$,
and 
\[  
\fra_{\vec{m}} \cdot
\fra_{\vec{\ell}}
\ \subseteq \ \fra_{\vec{m} +
\vec{\ell}} \]
for all $\vec{m} , \vec{\ell} \in \ZZ^r$. 
One puts $N^1(\fra\bull)_{\RR} = \ZZ^r
\otimes_{\ZZ}
\RR = \RR^r$, and
\[ \Nef(\fra\bull) \ \subseteq \
\Effbar(\fra\bull) \ \subseteq  \
N^1(\fra\bull)_{\RR} \] are defined as above.
Then $\text{Amp}(\fra\bull)$ and 
$\text{Big}(\fra \bull)$ are taken to be the
interiors of these cones. 

 Now fix a subvariety
$Z
\subseteq V$. Then under mild assumptions on
$\fra\bull$, one can mimic the global
constructions to define a continuous function
\[ \ord_Z : \text{Big}(\fra\bull) \lra \RR \]
that coincides with the global function when
$\fra\bull$ is the system of base-ideals just
discussed.\footnote{The assumption
is that $N^1(\fra \bull)$ should have a basis
consisting of 	``ample indices".} We refer to
\cite{AIBL}, \cite{Wolfe03} and
\cite{Wolfe05} for details.

 The
interesting point is that in this abstract
setting, the function $\ord_Z$
can be wild. In fact, Wolfe \cite{Wolfe03}
proves:
\begin{theorem}
There exist multi-graded families $\{
\fra\bull \}$ of monomial ideals on
$V = \CC^n$ for which the function
\[
\ord_0 : \textnormal{Big}(\fra\bull) \lra \RR 
\]
is nowhere differentiable on an open set. 
\end{theorem} 
\noi One conjectures of course that this sort
of behavior cannot occur in the global setting.
However all the known properties of $\ord_E$
can be established in this formal setting. 
So one is led to expect that there should be
global properties of these invariants that have
yet to be  discovered.

\def\cprime{$'$} \def\cprime{$'$} \def\cprime{$'$} \def\cprime{$'$}
  \def\cprime{$'$} \def\cprime{$'$} \def\cprime{$'$} \def\cprime{$'$}
  \def\cprime{$'$} \def\cprime{$'$} \def\cprime{$'$} \def\cprime{$'$}
\providecommand{\bysame}{\leavevmode\hbox to3em{\hrulefill}\thinspace}

\end{document}